\documentclass[11pt]{article}
\usepackage{subeqnarray,epsfig,amsfonts,amsthm,amssymb,comment}
\newcommand{\bs}{\begin{subeqnarray}}
\newcommand{\es}{\end{subeqnarray}}
\newcommand{\bss}{\begin{subeqnarray*}}
\newcommand{\ess}{\end{subeqnarray*}}
\newcommand{\ba}{\begin{eqnarray}}
\newcommand{\ea}{\end{eqnarray}}
\newcommand{\bas}{\begin{eqnarray*}}
\newcommand{\eas}{\end{eqnarray*}}
\newcommand{\bq}{\begin{equation}}
\newcommand{\eq}{\end{equation}}
\newcommand{\bqs}{\begin{equation*}}
\newcommand{\eqs}{\end{equation*}}
\newcommand{\la}{\label}
\newcommand{\sla}{\slabel}
\newcommand{\nn}{\nonumber}

\def\O{\Omega}

\def\p{\partial}

\def\na{\nabla}

\begin{document}
\title{The domain and property illusion of anomalous anisotropic electric conductivity}
\author{Kiwoon Kwon
\thanks{Department of Mathematics, Dongguk University-Seoul, 100715 Seoul, South Korea, 
         e-mail: kwkwon@dongguk.edu}}
\maketitle


\begin{abstract}
The unique determination of electrical conductivity is extensively studied for isotropic conductivity 
ever since Calderon's suggestion of the EIT (Electrical Impedance Tomography) problem.  
However, it is known that there are many anisotropic conductivities producing the same Dirichlet-to-Neumann map;
moreover the anisotropic conductivities giving the same Dirichlet-to-Neumann map are classified 
using the equivalence relation with respect to the change of variables.  
The change of variable argument is applied to the theory of near-cloaking: We are under an illusion that
the domain of anomaly is of a much smaller size than actually it is.
For this paper, we considered not only the illusion of the domain of the anomaly,
but also the illusion of the property of the anomaly when the background anisotropic conductivity is known.
\end{abstract}

\section{Introduction}
Recently, the idea of physical devices related to cloaking or invisibility is suggested more concretely than before
and draws much attentions from researchers in these areas \cite{ScienceLeonhardt, ScienceLeonhardtTyc, SciencePendry}. 
Invisibility and cloaking are closely related with the nonuniqueness of coefficients of the inverse problem modelled 
by EIT \cite{KangEIT, BryanLeise, GLU20031, GLU20032,  KohnEIT}, 
acoustic scattering \cite{KangIS,KohnIS, Nguyen}, 
electromagnetic scattering \cite{GKLU09, PRL},
and quantum scatteirng \cite{GKLU12,GKLU11}. 
In this paper, only the inverse problem modelled by EIT is considered.
 
Voltage potential in a Lipschitz domain $\O$ for an electrical current on the boundary $\p\O$ is given as follows:
\bs\la{eq:diffuse}
-\na\cdot(\sigma\na u) = 0  &&\qquad \mbox{ in } \O,
\sla{eq:diffuse1}\\
u = f  &&\qquad \mbox{ on } \p\O.
\sla{eq:diffuse2}
\es
If  $ f \in H^{1/2}(\p\O) $ and 
\bq\la{eq:bound}
0<L<\frac{y^t \sigma(x) y}{y^t y}< U 
\eq
for all $x\in \Omega,\; y\in\mathbb R^n\setminus\{0\}$ and some constants $L$ and $U$,
it is known to be $u\in H^1(\Omega)$ \cite{GilbargTrudinger}.
Therefore, we can define the Dirichlet-to-Neumann map $\Lambda_\sigma: H^{1/2}(\p\O) \rightarrow H^{-1/2}(\p\O)$ as $\Lambda_\sigma(f) = \nu\cdot (\sigma \nabla u)|_{\p\O}$ 
using the boundary trace operator, when $\Omega$ is a Lipschitz domain.

EIT is formulated to find $\sigma$  such that $\Lambda_\sigma=\Lambda$, given the Dirichlet-to-Neumann map $\Lambda$. In finite measurements case, EIT is to find  $\sigma$ such that $\Lambda_\sigma(f_i) = \Lambda(f_i), i=1,\cdots,N$ for given Dirichlet and Neumann boundary measurement pairs $(f_i,\Lambda(f_i))_{i=1,\cdots,N}$. 
 
Suppose that $D$ is a  Lipschitz domain compactly embedded in $\O$ such that $\O\setminus \overline D$  is connected. Reminding all the conductivities in Table 1 satisfies (\ref{eq:bound}), we classify conductivity $\sigma$ in (\ref{eq:diffuse}) into six cases in Table 2: Four cases of them (Case 1,2,4, and 5) includes the unknown obstacle $D$ inside the domain of interest $\O$, which have a different conductivity from the background $\O\setminus\overline D$. In other words, $\sigma$ is discontinuous along $\p D$ in these four cases.  

\begin{table}
\begin{tabular}{|c|l|c|l|}
\hline\\
$a$ & a known constant & $b$ & an unknown constant \\
\hline\\
$a(x)$ & a known function & $b(x)$ & an unknown function \\
\hline\\
$A$ & a known scalar matrix & $B$& an unknown scalar matrix \\
\hline\\
$A(x)$& a known function matrix & $B(x)$& an unknown function matrix\\
\hline
\end{tabular}
\caption{Electrical conductivities used in background and target obstacles} 
\end{table}

\begin{table}
\begin{tabular}{|c|c|c|c|}
\hline\\
 Case 1 & $\sigma = a + (b-a)\chi_D$ 
& Case 4 & $\sigma = A + (B-A)\chi_D$  
\\ \hline \\
Case 2 & $\sigma = a(x) + (b(x)-a(x))\chi_D$ 
&Case 5 & $\sigma = A(x) + (B(x)-A(x))\chi_D$\\
\hline \\
Case 3 & $\sigma = b(x)$
&Case 6 & $\sigma = B(x) $\\
\hline 
\end{tabular}
\caption{Electrical conductivities to be reconstructed on the whole domain of interest $\Omega$}
\end{table}

We explain briefly the results on the uniqueness and nonuniqueness in Section 2.
In Section 3, several results are given not only for illusion and cloaking of domain $D$ but also for property $B(x)$ illusion. 

\section{Uniqueness and nonuniqueness}
\subsection{Uniqueness}
Note that in the six cases in Table 2, the unknowns are $b,b(x),B, B(x)$ and $D$, whereas $a,a(x),A$ and $A(x)$ are assumed to be known. The uniqueness studies of EIT in these six cases proceeded in the following ways: 
\begin{itemize}
\item{Case 1 ($\sigma = a + (b-a)\chi_D$): \\
One or two measurements uniqueness results are known when $b$ is assumed to be known \cite{BFS, FI, KS}. The finite measurement uniqueness is known only for Case 1. In this case, the number of measurements, the geometry of the obstacle $D$, and the choice of suitable Dirichlet or Neumann data minimizing the number of measurements are interesting issues.}

\item{Case 2 ($\sigma = a(x) + (b(x)-a(x))\chi_D$) and Case 4  ($\sigma = A + (B-A)\chi_D$):\\
In \cite{IsakovEIT}, an orthogonality relation between the two solutions of (\ref{eq:diffuse}) for arbitary obstacles $D_1$ and $D_2$ is derived. Based on this relation and the Hahn-Banach theorem, a modified orthogonal relation is derived and the uniquness for Case 2 and Case 4 is proved. In Case 4, additional condition such that $B-A$ is positive-definite is used: this additional condition is generalized and removed by \cite{Ikehata} and \cite{KwonSheen}. The orthogonal and modified orthogonal equation is also used in Case 5.}

\item{Case 3 ($\sigma = b(x)$): \\
This is a generalization of Case 1 and Case 2, since we do not need to have the known background conductivity $a(x)$. Many mathematicians conduct extensive works in this case and thus, we are able to understand the unique determination of electrical conductivity when $b(x)\in L^\infty(\Omega)$ for the two-dimension $\cite{AstalaPaivarinta, BrownUhlmann, Nachmann}$ and $b(x)\in C^{3/2}(\Omega) $ for three and higher dimensions \cite{Alessandrini, Nachmann3D, SylvesterUhlmann87, SylvesterUhlmann88}.}

\item{Case 5 ($\sigma = A(x) + (B(x)-A(x))\chi_D$): \\
The uniqueness of $D$ is shown when $B(x)-A(x)$ is positive-definte \cite{IsakovEIT}; it is also presented when the positive-definite condition is removed in \cite{Kwon} and the nonuniqueness of $B(x)$ is pointed out by using the nonuniqueness result of Case 6 }

\item{Case 6 ($\sigma = B(x)$): \\
$B(x)$ is not unique for the given Dirichlet-to-Neumann map.  The number of $B(x)$ making the same Dirichlet-to-Neumann map depends on the number of the diffeomorphism on $\Omega$, which will be described in Definition 1. }
\end{itemize}

\subsection{Nonuniqueness in Case 6}
The nonuniqueness of EIT is observed early in \cite{KohnVogelius}. Before stating the result, let us define the diffeomorphism on $\O$ as follows:

{\bf Definition 1} A map $F$ defined almost everywhere from $\O$ into $\O$ itself is called diffeomorphism on $\O$, if the following conditions are satisfied:
\bs\la{eq:def}
&& F \mbox{ is bijective }\\
&& F \mbox{ is an identity on }\p\O\\
&& \mbox{Jacobian } DF \mbox{ is defined and bounded }\\
&& \det(DF) \neq 0
\es

{\bf Proposition 2} Let $F$ be a diffeomorphism on $\O$. Then, we have
\bq\la{eq:nonEIT}
\Lambda_\sigma = \Lambda_{F_*\sigma},
\eq
where
 \[ F_* \sigma = \frac {(DF) \sigma (DF)^t}{\det(DF)}.\]

{\bf Proof} 
Knowing the Dirichlet-to-Neumann map is equivalent to knowing the quadratic form  
\[ Q( v ) = \int \sigma \nabla v \cdot \nabla v \;\; v\in H^1(\O)\] 
by the polarization identity. 
Then, the proposition follows from the change of variable method for the quadratic form
\bq
   \int_\O \sum_{ij} \sigma_{ij} \frac{\p u}{\p x_i} \frac{\p u}{\p x_j} dx
 =\int_\O  \sum_{ij} \sigma_{ij} \frac{\p u}{\p y_k} \frac{\p y_k}{\p x_i}
                                                   \frac{\p u}{\p y_l}  \frac{\p y_l}{\p x_j} 
                                                   \det\left(\frac{\p x}{\p y}\right) dy  
\eq 
as explained in \cite{KohnEIT}. \qed

$F_*$ in Proposition 2 is called the 'push forward map'  for a diffeomorphism $F$ on $\O$. After Proposition 1 is known, the following questions in Case 5 and Case 6 are raised: 
\begin{itemize}
\item{ Is the change of variable (\ref{eq:nonEIT}) a unique obstruction to the uniqueness for Case 6? 
In other words, does $\Lambda_{\sigma_1} = \Lambda_{\sigma_2}$ imply $\sigma_2 = F_* \sigma_1$ for some diffeomorphism $F$ on $\O$? }
\end{itemize}

The answer is yes in Case 6, which means that there are no nonuniqueness cases other than the change of variables. Many researchers, including Gunter Uhlmann, proved this problem for the following subcases:
\begin{itemize}
\item{ $n=2$ and $\sigma\in C^3(\O)$ \cite{Sylvester}}
\item{ $n=2$ and $\sigma\in C^{0,1}(\O)$(Lipschitz functions) \cite{SunUhlmann}}
\item{ $n=2$ and $\sigma\in L^\infty(\O)$ \cite{AstalaPaivarintaLassas}}
\item {$ n\ge 3$ : $\p D$ and $B$ are analytic \cite{LassasUhlmann, LeeUhlmann} }
\end{itemize}

The answer for Case 5 is given in the following section.

\section{Illusions of domain ($D$) and property ($B$) in Case 5 and Case 6}

In Case 5 ($\sigma = A(x) + (B(x)-A(x)) \chi_D$), the uniqueness of $D$ is solved in \cite{Ikehata, IsakovEIT, Kwon}; yet, the nonuniqueness of $B(x)$
inside $D$ is only commented in \cite{Kwon} using the result for Case 6 (Section 2.2). Theorem 1 in  \cite{Kwon} is rephrased in the following theroem:

{\bf Theorem 3}
Assume that  
\bq\la{eq:th3a}
 A(x), A^{1/2}(x), B_i(x), B_i^{1/2}(x) \mbox{ are Lipschitz continuous on }\O 
\eq
and there is a jump on $\p D$ such that 
\bq\la{eq:th3b}
\det(B_i) \neq \det(A)  \mbox{ for } n=2 \mbox{ and },
B_i\neq A \mbox{ for } n\ge 3 
\eq
for $i=1,2$.
Then, $\Lambda_{\sigma_1} = \Lambda_{\sigma_2}$ implies 
\bq\la{eq:th3c1}
 D_1 = D_2 .
\eq
Denoting $D=D_1= D_2$,  on the boundary $\p D$, we have 
\bq\la{eq:th3c2}
\det(B_1) = \det(B_2) \mbox{  for } n=2, \mbox{ and }
         B_1    =  B_2       \mbox{  for } n\ge 3.
\eq
If we further assume that
\ba\la{eq:th3d}
 &&B_1 = B_2  \mbox{ on }\p D \mbox{ for } n=2, 
\sla{eq:th3d1} \\
 && \p D \mbox{ is connected and analytic and } B|_D \mbox{ is analytic  for } n\ge 3. 
\nn
\ea
Then, there is a diffeomorphism $F$ on $D$ such that  
\bq\la{eq:th3e}
    B_2 = F_* B_1. 
\eq

From Theorem 3, in Case 5, the domain $D$ is uniquely determined by the Dirchlet-to-Neumann map; however, the property $B$ is nonunique for the same Dirichlet-to-Neumann map. 
Specifically, the domain of anomaly cannot be cloaked into a much smller domain, but
the property could be illued into other property. In other words, in Case 5, the size illusion is impossible 
but the property illusion is possible.

\subsection{Illusion of domain $D$ for Case 6}
Consider the conductivity 
      \[\sigma = A(x) + (B(x)-A(x))\chi_D.\]
In Case 5, $A(x)$ is the known property. Then, by Theorem 3, there is no cloaking or illusion for domain $D$.  However, if $A(x)$ is assumed to be unknown, $\sigma$ is a special case of Case 6 with $\sigma \in L^\infty(\O)$. In this case, there is illusion or cloaking of domain  $D$. Specifically, there is a diffeomorphism $F$ on $\O$ which illude the domain $D$ into arbitrary small domain $U_\epsilon$ through the following conductivity:  
\[
    \sigma_F(x) = 1 + (F_* B (x)-1) \chi_{U_\epsilon}(x)
\]
For example, if $\O,D$ and $U_\epsilon$ are balls of radius 2, 1, and $\epsilon>0$, respectively, centered at $\mathbb O$, we could take the diffeomorphism $F_\epsilon$ on $\O$ as the inverse of the followings :
 \bq\la{eq:F}
 F_\epsilon(x) = \epsilon x \chi_D (x) + \left( (2-\epsilon)|x| - (2-2\epsilon)\right) \frac{x}{|x|}  \chi_{\O\setminus \overline D}(x) \eq
where $F_\epsilon$ maps $D$ into $U_\epsilon$. If $\epsilon$ goes to $0$, we stae that $D$ is "perfectly cloaked" by the corresponding diffeomorphism $F_0$.

Let us interpret the cloaking in a different point of view. Let $F$ be a diffeomorphism on $\O$ such that $F(D_1) = D_2$. Define three kinds of conductivities as follows:
\bs
\sigma_1 &=& I_n\chi_{\O\setminus\overline{D_1}} + B_1\chi_{D_1},  \slabel{eq:1}\\
\sigma_2 &=& I_n\chi_{\O\setminus\overline{D_2}} + B_2\chi_{D_2}. \slabel{eq:2}\\
\sigma_3 &=& A_3 \chi_{\O\setminus\overline{D_3}} + B_3\chi_{D_3},  \slabel{eq:3} 
\es
where $I_n$ is a $n-$ dimensional identity matrix.

{\bf Theorem 4} Suppose that 
\bq\la{eq:th3_1}
 D_1 \neq D_2 
\eq
and (\ref{eq:th3a}), (\ref{eq:th3b}) hold.
Then
\bq\la{eq:th3_2}
 \Lambda_{\sigma_1} \neq \Lambda_{\sigma_2} 
\eq
Suppose that
\bq\la{eq:th3_3}
 A_3 = F_* 1,\;  B_3 = F_* B_1, \;{\rm and }\; D_3 = F_* D_1 
\eq
for some diffeomorphism $F$ on $\O$.
Then
\bq\la{eq:th3_4}
\Lambda_{\sigma_1} = \Lambda_{\sigma_3} 
\eq
Conversely, if (\ref{eq:th3_4}) holds in the two dimension, then there is diffeomorphism $F$ on $\O$ 
satisfying (\ref{eq:th3_3}).

{\bf Proof}
(\ref{eq:th3_2}) and (\ref{eq:th3_4}) are followed by Theorem 3 and Proposition 2.
The converse is followed by \cite{AstalaPaivarintaLassas}. \qed

From (\ref{eq:th3_2}), we do not have any illusion of domain $D_1$ to another domain $D_2$ in Case 5, when the background conductivity $A$ is known. However, if we do not know the background conducitvity $A$ in Case 6, we could illude domain $D_1$ into a more smaller domain $D_3$  ( or perfectly cloak $D_1$ ) by the diffeomorphism $F$ on $\O$ , for example $F_\epsilon$ in (\ref{eq:F}),  using  (\ref{eq:th3_4}) in Theorem 4.

\subsection{Illusion of property $B$ for Case 5}
Although domain ($D$)  illusion is impossible in Case 5, property ($B$) illusion is possible
by the following theorem:

{\bf Theorem 5}
Let 
\bq\la{eq:th4_1}
  \sigma_i = I_n + (B_i - I_n) \chi_D, \; i=1,2.
\eq
If 
\bq\la{eq:th4_2}
  B_2 = F_* B_1
\eq
for some diffeomorphism $F$ on $\O$.
Then, 
\bq\la{eq:th4_3}
\Lambda \sigma_1 = \Lambda \sigma_2.
\eq
Conversely, if (\ref{eq:th4_3}) holds, then there is diffeomorphism $F$ on $D$ such that (\ref{eq:th4_2}) holds under the conditions (\ref{eq:th3a}),(\ref{eq:th3b}).and (\ref{eq:th3d}). 

{\bf Proof}
 (\ref{eq:th4_3}) follows Proposition 2.
And the existence of the diffeomorphism $F$ on $D$ follows (\ref{eq:th3e}) in Theorem 3. \qed

From Theorem 5, property illusion is possible up to the change of variables by diffeomorphisms on $\O$. 
In \cite{PRL}, a few interesting property illusions are considered; The optical transformation of an object into another object is 
considered. In the paper, stereoscopic image of a man is transformed into an illusion image of a woman and dielectic spoon of the electric permeability 2 into an illusion image of metallic cup of electric permeability -1 in an electromagnetic scattering problem.

\section{Conclusions}
In EIT, the uniqueness for isotropic conductivities is summarized. It is shown that the 
nonuniqueness for anisotropic conductivities is possible if and only if the conductivity is pulled back
by the diffeomorphism on the domain of interest. 
In Case 6, when the background is also unknown, illusion and perfect cloaking of domain is possible.  
In Case 5, when the background is known, illusion of domain is impossible; yet,the illusion of property inside the domain is possible by the pull-back of the diffeomorphism on the domain.  

\section{Acknowledgment}
This work was supported by Dongguk University Research Fund of 2010.

\bibliography{rem1}
\bibliographystyle{plain}

\end{document}